\documentclass{gtmon_a}
\pdfoutput=1

\usepackage[T1,T5]{fontenc}

\proceedingstitle{Proceedings of the School and Conference in Algebraic
Topology (The Vietnam National University, Hanoi, 9-20 August 2004)}
\conferencestart{9 August 2004}
\conferenceend{20 August 2004}
\conferencename{School and Conference in Algebraic Topology}
\conferencelocation{Vietnam National University, Hanoi, Vietnam}

\editor{John Hubbuck}
\givenname{John}
\surname{Hubbuck}

\editor{Nguy\~\ecircumflex{}n H V H\uhorn{}ng}
\givenname{H\uhorn{}ng}
\surname{Nguy\~\ecircumflex{}n}

\editor{Lionel Schwartz}
\givenname{Lionel}
\surname{Schwartz}

\title[Characteristic classes for extraspecial $p$--groups] 
{Evens norm, transfers and characteristic classes\\for extraspecial $p$--groups}
\author{Ph\d{a}m Anh Minh}
\surname{Ph\d{a}m}
\givenname{Minh}

\address{Department of Mathematics\\
College of Science\\
University of Hue\\\newline
Dai hoc Khoa hoc\\
77 Nguyen Hue\\
Hue\\
Vietnam}
\email{}
\urladdr{}

\volumenumber{11}
\issuenumber{}
\publicationyear{2007}
\papernumber{9}
\startpage{179}
\endpage{200}

\doi{}
\MR{}
\Zbl{}

\arxivreference{}



\received{30 November 2004}
\revised{}
\accepted{23 January 2005}
\published{14 November 2007}
\publishedonline{14 November 2007}
\proposed{}
\seconded{}
\corresponding{}
\editor{}
\version{}

\renewcommand \inf{\mathrm{Inf}}
\renewcommand \Re{\mathbb{R}}
\renewcommand \P{\mathbb{P}}

\def\A{{\mathcal A}}
\def\F{{\mathbf F}}
\def\N{{\mathcal N}}
\def\T{{\mathcal T}}
\def\U{{\cal U}}

\newcommand{\im}{\mathop{\rm Im}}
\newcommand{\ra}{\longrightarrow}
\newcommand{\sra}{\rightarrow}

\def\ka{\kappa}
\def\kapi{\kappa_{n,r}^\iota}
\def\res{{\rm Res}}
\def\tr{{\rm tr}}
\def\KG{{K \rightarrow G}}

\newcommand{\norm}[2]{\mathcal{N}_{#1\sra #2}}

\renewcommand{\theenumi}{\roman{enumi}}

\makeatletter
\def\cnewtheorem#1[#2]#3{\newtheorem{#1}{#3}[section]
\expandafter\let\csname c@#1\endcsname\c@thm}

\theoremstyle{plain}
\newtheorem{thm}{Theorem}[section]
\cnewtheorem{prop}[thm]{Proposition}
\cnewtheorem{cor}[thm]{Corollary}
\newtheorem{lem}[thm]{Lemma}

\theoremstyle{definition}
\cnewtheorem{defn}[thm]{Definition}
\cnewtheorem{conj}[thm]{Conjecture}
\cnewtheorem{exmp}[thm]{Example}

\theoremstyle{remark}
\newtheorem*{rem}{Remark}

\makeatother  %  move after \newtheorem block

\begin{document}

\begin{htmlabstract}
Let P be the extraspecial p&ndash;group of order p<sup>2n+1</sup>,
of p&ndash;rank n+1, and of exponent p if p>2. Let Z be the center of
P and let &kappa;<sub>n,r</sub> be the characteristic classes of degree
2<sup>n</sup>-2<sup>r</sup> (resp. 2(p<sup>n</sup>-p<sup>r</sup>))
for p=2 (resp. p>2), 0&le;r&le;n-1, of a degree p<sup>n</sup> faithful
irreducible representation of P. It is known that, modulo nilradical,
the &iota;th powers of the &kappa;<sub>n,r</sub>s belong to
T=Im(H<sup>*</sup>(P/Z,<b>F</b><sub>p</sub>)/&radic;0&rarr;H<sup>*</sup>(P,<b>F</b><sub>p</sub>)/&radic;0),
with &iota;=1 if p=2, &iota;=p if p>2.  We obtain formulae
in H<sup>*</sup>(P,<b>F</b><sub>p</sub>)/&radic;0 relating
the &kappa;<sub>n,r</sub><sup>&iota;</sup> terms to the
ones of fewer variables. For p>2 and for a given sequence
r<sub>0</sub>,&hellip;,r<sub>n-1</sub> of non-negative
integers, we also prove that,  modulo-nilradical, the element
&prod;<sub>r<sub>i</sub></sub>&kappa;<sup>r<sub>i</sub></sup><sub>n,i</sub>
belongs to T if and only if either r<sub>0</sub>&ge;2,
or all the r<sub>i</sub> are multiple of p. This gives the
determination of the subring of invariants of the symplectic group
Sp<sub>2n</sub>(<b>F</b><sub>p</sub>) in T.
\end{htmlabstract}

\begin{abstract} 
Let P be the extraspecial $p$--group of order $p^{2n+1}$, of $p$--rank
$n+1$, and of exponent $p$ if $p>2$. Let $Z$ be the center of $P$ and
let $\kappa_{n,r}$ be the characteristic classes of degree $2^n - 2^r$
(resp.\ $2(p^n-p^r)$) for $p=2$ (resp.\ $p>2$), $0\leq r \leq n-1$, of a
degree $p^n$ faithful irreducible representation of $P$. It is known that,
modulo nilradical, the $\iota$th powers of the $\kappa_{n,r}$'s belong to
${\mathcal T}={\rm Im}\bigl(H^*(P/Z,\mathbf{F}_p)/\sqrt 0 \smash{\buildrel{{\rm
Inf}} \over{\rightarrow}} H^*(P,\mathbf{F}_p)/\sqrt 0\bigr)$, with $\iota= 1$ if
$p=2$, $\iota= p$ if $p>2$.  We obtain formulae in $H^*(P,\mathbf{F}_p)/\sqrt
0$ relating the $\kappa_{n,r}^\iota$ terms to the ones of fewer
variables. For $p>2$ and for a given sequence $r_0,\ldots,r_{n-1}$
of non-negative integers, we also prove that,  modulo-nilradical,
the element $\prod_{r_i}\kappa^{r_i}_{n,i}$ belongs to ${\mathcal T}$
if and only if either $r_0 \geq 2$, or all the $r_i$ are multiple of
$p$. This gives the determination of the subring of invariants of the
symplectic group $Sp_{2n}(\mathbf{F}_p)$ in ${\mathcal T}$.
\end{abstract}

\maketitle

\section[Introduction]{Introduction}
%\label{sec:sec1}

Let $p$ be a prime number. For a given group $P$ , denote by
$H^*(P )$ the mod--$p$ cohomology algebra of $P$. We are 
interested
in the case where $P = P_n$, the extraspecial group $p$--group of
order $p^{2n+1}$, of $p$--rank $n+1$, and of exponent $p$ if 
$p>2$.
It is known (see the work of Green--Leary \cite{GL} 
or Quillen\cite{Q}) that, for $p>2$
(resp.\ $p=2$), there are exactly $n$ Chern (resp.\ Stiefel--Wihitney)
classes $\kappa_{n,r}$ of degree $2(p^n-p^r)$
(resp.\ $2^n - 2^r$), $0\leq r \leq n - 1$, of a degree $p^n$
faithful irreducible representation of $P$; these classes
restrict to maximal elementary abelian subgroups of P as Dickson
invariants.

Set $E = E_n = P/Z$, with $Z$ the center $P$, $E$ is then a 
vector
space of dimension $2n$ over $\F_p$. Set $h^*(P ) = H^*(P )/\sqrt
0$  (so $h^*(P ) = H^*(P )$ for $p=2$, by \cite{Q} and denote by
${\mathcal T} = {\mathcal T}_n$ the subring of $h^*(P)$ equal to
the image of the inflation ${\rm Inf}^E_P$, modulo nilradical. 
For
$p=2$, it follows from \cite{Q} that all the $\kappa_{n,r}$ terms
belong to ${\mathcal T}$. For $p>2$, this fact does not hold, as
shown by Green and Minh \cite{G,GM}; however, in \cite{GM}, it is 
also proved that all $p$th powers of the  
$\kappa_{n,r}$ terms, are in
turn, belonging to ${\mathcal T}$.

For convenience, set  $\iota= 1$ for $p=2$, and $\iota = p$ for
$p>2$. It follows that, for $0 \leq  r \leq n - 1$, there exists
$f_{n,r} \in H^*(E)$ such that  ${\rm Inf}(f_{n,r})\doteqdot
\kapi$. Here and in what follows $a \doteqdot b$ means $a=b$
modulo $\sqrt 0$. Via the inflation map, as elements of $h^*(P)$,
$f_{n,r}$ can be identified with $\kapi$, $0 \leq r \leq n-1$. 
The
first aim of this paper is to get an alternating formula
expressing $f_{n,r}$ by means of $f_{n-1,r-1}$ and $f_{n-1,r}$.
This work is motivated from the elegant formula, for $p>2$
(resp.\ $p=2$), expressing Chern (resp.\ Stiefel--Whitney) classes of the
regular representation $r_A$ of the elementary abelian $p$--group
$A$ via such classes of fewer variables. It is known that, if $A$
is of rank $m$ and $p>2$ (resp.\ $p=2$), then the $2(p^ m -
p^r)$th Chern (resp.\ $(2^m - 2^r)$th Stiefel--Whitney) class
of $r_A$ is the Dickson invariant $Q_{m,r}$ of the same degree
with variables  in a basis $ x_1, \ldots , x_m$ of $\beta H^1(A)$
(resp.\ $ H^1(A))$, with  the  Bockstein homomorphism. These
invariants are related by
$$
Q_{m,r} = Q_{m-1,r} V_m^{p-1}+Q^p_{m-1,r-1},
$$
where 
$$ 
V_m = \prod_{\lambda_i \in \F_p}
\left(\lambda_1x_1+\cdots+\lambda_{m-1}x_{m-1}+x_m\right)
=
(-1)^{m-1}\sum_{s=0}^{m-1}(-1)^sQ_{m-1,s}x_m^{p^s},
$$ 
is the M\`ui invariant.

In so doing, we need to use the Evens norm and
transfers from maximal subgroups of $P$. Some basic properties of
the Evens norm, in the relation with modular invariants, are
recalled in \fullref{sec:sec2}. In \fullref{sec:sec3}, 
we show how to obtain characteristic classes of $P$ by means 
of the Evens norm (\fullref{thm:thm3.7}). 
\fullref{thm:thm3.8} describes the image of such classes via the
Evens norm. From this, we obtain formulae relating characteristic
classes with such classes of fewer variables (\fullref{cor:cor3.9}).

Let $r_0 , \ldots , r_{n-1}$ be a sequence of non-negative 
integers.
In \fullref{sec:sec4}, we prove that, for $p>2$, modulo 
nilradical, the product $\prod_{i\geq 0} \kappa_{n,i}^{r_i}$ 
belongs to ${\mathcal T}$ if and only if  
either $r_0 \geq 2$, or all the $r_i$ are
multiple of $p$ \fullref{thm:thm4.1}. This generalizes a result, given 
by Green and Leary \cite{G,GL}, 
proving that $\kappa^{s}_{n,0}$ belongs to
${\mathcal T}$ provided either $s \geq 2^n$, or $s\geq 2$ and 
$n \leq 2$. As a consequence, we obtain in the last section the
determination of the subring of invariants of the symplectic 
group in ${\mathcal T}$ \fullref{thm:thm5.1}.

For convenience, given a subgroup $K$ of a group 
$G$, any element of $H^*(G)$ is also considered as an element of
$H^*(K)$ via the restriction map ${\rm Res}^G _K$. Also, if $K$ 
is normal in $G$, then any element of $H^*(G/K)$ can be considered
as an element of $H^*(G)$ via the inflation map 
${\rm Inf}^{G/K}_K$.

\section{Evens norm and M\`ui invariants}
\label{sec:sec2}

Given a polynomial algebra $F = \F_p[t_1, \dots , t_k]$ and 
$1\leq m \leq  k - 1$, define the M\`ui invariant \cite{Mu} 
\begin{equation}
V_{m +1} = V_{m+1}(t_1, \ldots , t_{m +1} ) =
\prod_{\lambda_i \in \F_p}
(\lambda_1t_1+\cdots+\lambda_{m}t_{m}+t_{m+1}).
\label{eqn:eqn1}
\end{equation}

It follows from the work of L\,E Dickson \cite{D} that 
$$ 
V_{m +1}=
(-1)^{m}\sum_{s=0}^{m}(-1)^sQ_{m,s}t_{m+1}^{p^s}
$$
with
$Q_{m,s} =Q_{m,s}(t_1, \ldots , t_{m})$ the Dickson invariants defined
inductively as follows (we shall omit the variables, if no
confusion can arise).
\begin{align*}
Q_{m,m} & = 1 \\
Q_{m,0} & = \prod_{\substack{\lambda_j \in \F_p\\
\text{$\lambda_j$ not all equal $0$}}}
(\lambda_1 t_1 + \cdots + \lambda_m t_m) \\
Q_{m,s} & = Q_{m-1,s} V_m^{p-1} +
Q^p_{m-1 ,s-1}.
\end{align*}
By \eqref{eqn:eqn1} the $Q_{m,s}$ are independent of the choice of
generators $t_1,\ldots,t_m$ of $\F_p[t_ 1,\ldots,t_m]$.
Hence, if $(t_1, \ldots , t_m)$ is a basis of $H^1(W )$ 
(resp.\ $\beta H^1(W )$) with $W$ an elementary 
abelian $2$--group (resp.\ $p$--group with $p > 2$) 
of rank $m$, we may write 
\begin{align*}
Q_{m,s}(t_1,\ldots, t_m) & = Q_s(W ) \\
V_{m+1}(t_1, \ldots, t_k, X) & = V (W,X).
\end{align*}
The M\`ui invariants can be obtained by means of Evens norm map
$\norm UW$
with $U$ a subgroup of $W$ (see \fullref{cor:cor2.2} below). Let us
recall that, for every maximal subgroup $K$ of a $p$--group $G$,
and for $\xi \in H^r(K)$, we may define the Evens norm map
$$
\norm KG (\xi) \in  H^{pr}(G).
$$
Here are some properties of $\norm KG$. For details of the proof,
the reader can refer to the work of Evens \cite{E},
Minh \cite{M2} or M\`ui \cite{Mu}.

\begin{prop} 
\label{prop:prop2.1}
Let $G$, $G'$ be $p$--groups and let $K$ be a 
subgroup
of $G$.
\begin{enumerate}
\item
If $N$ is a subgroup of $K$, then 
$$
{\mathcal N}_{N\rightarrow G} = {\mathcal N}_{K\rightarrow G}\circ
{\mathcal N}_{N\rightarrow K}.
$$
\item
If $H$ is a subgroup of $G$ and $G=\coprod_{x \in D} KxH$, then, 
for $\xi \in H^r(G)$,
$$
\res^G_H {\mathcal N}_{K\rightarrow G}(\xi)=
\prod_{x \in D}{\mathcal N}_{H \cap {}^xK\rightarrow H}
\res^{{}^xK}_{H\cap {}^xK}({}^x \xi).
$$
\item
If $K$ is a subgroup of $G'$ and $f \co G' \rightarrow G$
is a homomorphism such that $f(K' ) \subset K$ and $f$
induces a bijection $G'/K' \cong G/K$ of coset spaces, then, for
$\xi \in H^r(K)$,
$$
{\mathcal N}_{K'\rightarrow G'} (f_{|K'})^*(\xi)=f^* 
\left({\mathcal N}_{K\rightarrow G}(\xi)\right).
$$
In particular, if $N$ is a normal subgroup of $G$ and $N \subset
K$, then, for $\xi \in H^r(K/N)$, 
$$
{\mathcal N}_{K\rightarrow G} \inf_K^{K/N}(\xi)= 
\inf_G^{G/N}{\mathcal N}_{K/N\rightarrow G/N}(\xi).
$$
\item
If $\xi, \xi' \in  H^n(K)$, then 
$$
\N_{K\rightarrow G}(\xi +\xi' ) = 
\N_{K\rightarrow G}(\xi  )+\N_{K\rightarrow G}( \xi' )
$$
modulo a sum of transfers from proper subgroups of $G$ containing
the intersection of the conjugates of $K$.
Hence, the norm map is
in general non-additive, although $\N_\KG \circ \res^G_K$ is.
\item
If   $\xi \in H^r(K)$, $\xi' \in  H^s(K)$, and $[G : K] = n$,
then
$$
\N_\KG(\xi.\xi') = 
(-1)^{\frac{n(n-1)}{2} rs} \N_\KG(\xi) \N_\KG(\xi' ).
$$
\item
Assume that $G = K \times E$, with $E=(\F_p)^m$. Consider $E$
as the group of all translations on a vector space $S$ of
dimension $m$ over $\F_p$ and let
$W (m)$ be an $E$--free acyclic complex with augmentation 
$\epsilon \co W(m) \rightarrow \F_p$. Let $C$ be a cochain complex of
which the cohomology is $H^*(K)$ and set 
$C^S = \otimes_{c \in S} C_c$, $C_ c = C$. 
Then 
$$
\N_{K \rightarrow K \times E} = d_m^*P_m,
$$
where 
$P_m \co H^r(C)\rightarrow H^{p^m r}_E\left(W(m) \otimes C^S \right)$ 
is the Steenrod power map, and 
$$d^*_m \co H^{p^mr}_E \left(W(m) \otimes C^S\right) 
\rightarrow H^*(E) \otimes H^*(C)$$is induced by the diagonal $C \rightarrow C^S$ and the K\"unneth
formula.
\end{enumerate}
\end{prop}

In the rest of this section, suppose that $W$ is an elementary
abelian $p$--group of rank $n + 1$ and $U$ a subgroup of $W$ of
index $p^m$.

By \fullref{prop:prop2.1}(vi), 
$\N_{U \rightarrow W}=d^*_mP_m$. 
The first part of the following corollary is then
originally due to M\`ui \cite{Mu} and reproved by Okuyama and Sasaki
\cite{OS}; the second one was given by H\uhorn{}ng and Minh 
\cite[Proof of Theorem B]{HM}.

\begin{cor}
\label{cor:cor2.2}
For $p = 2$ (resp.\ $p > 2$) and for every $x \in 
H^1(W)$ (resp.\ $x \in \beta H^1(W)$),
\begin{enumerate}
\item
$\N_{U \rightarrow W}\bigl({\rm Res}^W_U(x)\bigl)= V (W/U, x)$
\item
with $T$ the maximal subgroup of $W$ satisfying 
$\res^W_T (x) = 0$ then
\begin{multline*}
(-1)^r\norm TW (Q_r(T )) =
\sum_{i=r}^{n}(-1)^iQ^p_i (T)x^{p^{i+1}-p^{r+1}} \\
- \biggl[\sum_{i=r+1}^{n}(-1)^iQ_i(T)x^{p^i-p^{r+1}} \biggr]
\biggl[\sum_{i=0}^{n}(-1)^iQ_i(T )x^{p^i}\biggr]^{p-1}.
\end{multline*}
\end{enumerate}
\end{cor}

In the following corollary, $G$ is supposed to be a $p$--group
given by a central extension and 
$$ 
\{0\} \ra \Z/p \ra G \buildrel{ {j }}\over {\ra} W \ra \{0\} 
$$
and $K=j^{-1}(U)$. Set
$$
H^{ev}(U) =
\begin{cases}
H^*(U) & p=2, \\
\sum_{n\geq 0} H^{2n}(U) & p>2.
\end{cases}
$$

The following is straightforward from \fullref{prop:prop2.1}.
\begin{cor} 
\label{cor:cor2.3}
The composition map
$$
H^{ev}(U) \buildrel{ {{\rm Inf}^U_K }}\over {\ra}H^*(K) 
\buildrel{{\N_\KG }}\over {\ra}H^*(G)
$$
is a ring homomorphism. 
\end{cor}

In \cite{M2} we proved the following proposition.
\begin{prop}
\label{prop:prop2.4}
Let $\xi \in H^q(G)$. Set $\mu(q)=(-1)^{hq}h!$ with $h=(p-1)/2$
for $p>2$. If $K=\ker (u)$ with $u \in H^1(G)$, $ u \not = 0$, 
then, by setting
$v=\beta(u)$, we have
$$   
\norm KG (\res^G_K(\xi)) = \begin{cases}
\sum_i Sq^i(\xi)u^{q-i} & $p=2$ \\
\mu(q) \displaystyle\sum_{\substack{\epsilon = 0,1 \\
0 \leq 2i \leq q- \epsilon }}
(-1)^{\epsilon+i}\beta^\epsilon {\mathcal P}^i
(\xi)v^{(q-2i)h-\epsilon}u^\epsilon & $p>2$. 
\end{cases}
$$
where $Sq^i$ (resp.\ ${\mathcal P}^i$) denotes the Steenrod
operation for $p=2$ (resp.\ $p>2$).
\end{prop}

\section[Characteristic classes]
{Characteristic classes for extraspecial $p$--groups }
\label{sec:sec3}

Let $E=E_n$, $n\geq 1$, be the elementary abelian $p$--group of
rank $2n$. Let $x_1,\ldots,x_{2n}$ be a basis of 
$H^1(E) = {\rm Hom}(E, \F_p)$ and define
$$   
y_i = 
\begin{cases}
x_i & p=2 \\
\beta(x_i) & p>2.
\end{cases}
$$
$1 \leq i \leq  2n$, with $\beta$ the Bockstein homomorphism. We
have
$$   
H^*(E) = 
\begin{cases}
\F_p[y_1,\ldots,y_{2n}] & p=2 \\
\Lambda[x_1,\ldots,x_{2n}] \otimes \F_p[y_1,\ldots,y_{2n}] & p>2,
\end{cases}
$$
with $\Lambda[s, t,\ldots ]$ (resp.\ $\F_p[s, t,\ldots ]$) the
exterior (resp.\ polynomial) algebra with generators 
$s, t, \ldots$
over $\F_p$. Let $P = P_ n$ be the extraspecial $p$--group given 
by the central extension 
$$
\{1\} \ra Z/p \buildrel{ {i }}\over {\ra}P_n \ra E \ra \{1\}
$$
classified by the cohomology class 
$x_ 1x_2 + \cdots + x_{2n-1}x_{2n} \in H^2(E)$. 
The following notation will be used.
Set $Z=i(\Z/p)$, the center of $P$. For every elementary subgroup
$A$ of $P$ containing $Z$, write $A/Z=Z'$, so $A=A' \times Z$, 
and $A'$ is of rank $n$ if $A$ maximal elementary abelian in $P$. Fix
a generator $\gamma$ of $H^1(Z)$  (resp.\ $\beta H^1(Z)$) for 
$p = 2$ (resp.\ $p > 2$). This element, and also every element of 
$H^*(A' )$, are then considered as
elements of $H^*(A)$ via the inflation maps.

Denote by ${\mathcal A}$ the set of maximal elementary abelian
subgroups of $P$. Set $h^*(P)=H^*(P) / \sqrt 0$. By the work of Quillen
\cite{Q}, the map induced by the restrictions
$$
h^*(P)\buildrel{ {\res }}\over {\ra} \prod_{A \in {\mathcal A}}
h^*(A) 
$$ 
is injective. Therefore the maps
$$
h^*(E) \buildrel{ {\inf }^E_P}\over {\ra} h^*(P )
\qquad\text{and}\qquad
h^*(E)\buildrel{ {\res }}\over {\ra}\prod_{A \in \mathcal A} 
h^*(A' )
$$
have the same kernel. Let ${\mathcal T}={\mathcal T}_n$ be the
subring of $h^*(P )$ equal to the image of the inflation $\inf_P^E$.
For elements $\xi$, $\eta$ of $h^*(E)$, it follows that 
$\inf^E_P(\xi) = \inf^E_P(\eta)$ if and only if 
$\res^E_{A'}(\xi) = \res^E_{A'}P(\eta)$, for every $A \in \mathcal A$.

We are now interested in Chern (resp.\ Stiefel--Whitney) classes,
for $p> 2$ (resp $p=2$), of a degree $p^n$ faithful irreducible
representation of $P$. Fix a nontrivial linear character $\chi$ 
of $Z$. We have then an irreducible character $\hat \chi$ of $P$
given by
$$   
\hat \chi (g) = 
\begin{cases}
p^n \chi(g) & g \in Z \\
0 & {\rm otherwise}.
\end{cases}
$$
Let $\rho$ be a representation affording the character $\hat\chi$ and set
\begin{eqnarray*}
\zeta=\zeta_n & = &
\begin{cases}
c_{p^n}(\rho)& p>2\\
w_{2^n}(\rho) & p=2,
\end{cases} \\
\kappa_{n,r} & = & 
\begin{cases}
 (-1)^{n-r}c_{p^n-p^r}(\rho) & p>2 \\
w_{2^n-2^r}(\rho) & p=2,
\end{cases}
\end{eqnarray*}
$0 \leq r \leq n$. We have the following theorem.

\begin{thm}[Green--Leary \cite{GL}, Quillen \cite{Q}]\quad
\label{thm:thm3.1}
\begin{enumerate}
\item
In $h^*(P)$,
$$  1-\kappa_{n,n-1}+\cdots +(-1)^{n}\kappa_{n,0}+\zeta_n= 
\begin{cases}
c(\rho) &  $p>2$ \\
w(\rho) & $p=2$,
\end{cases}
$$
the subring of $h^*(P)$ generated by non-nilpotent Chern classes
is generated by 
$$
\kappa_{n,0},\ldots,\kappa_{n,n-1},y_1,\ldots,y_{2n},\zeta_n. 
$$
\item
For every $0 \leq i \leq n-1$ and for every $A \in \mathcal A$
\begin{align*}
\res^P_A(\zeta)& =V(A',\gamma) \\
\res^P_A(\kappa_{n,i}) & = Q_i(A').
\end{align*}
\end{enumerate}
\end{thm}

In the article \cite{GM} by Green and Minh, Chern classes of $P$ are
also obtained by means of the inflation $\inf^E_P$ and
transfer maps $\tr^K_P$ with $K$ maximal in $P$. Similar
results for the case $p = 2$ can also be obtained by using the
same argument. The result can be stated as follows. Let $x$ be a
non-zero element of $H^1(P )$ and set $H_x = \ker (x)$. Pick a
rank one subgroup $U \not = Z$ of the center of $H_x$.  So $H_x=
P_{n-1} \times U$. By the K\"unneth formula, we can consider any
element of $H^*(P_{n-1})$ (and of $H^*(U)$) as an element of
$H^*(H_x)$. For $0 \leq r \leq n-1$, set
$$  
\chi_{r,x}=   
\begin{cases}
\tr_P^{H_x}(\kappa_{n-1,r}\zeta^{p-1}_{n-1}) & n\geq 2\\
\tr_P^{H_x}(\zeta_{n -1}^{p-1}) & n=1.
\end{cases}
$$
For $r\geq 0$, define
$$  
z^{(r)}_n=   
\begin{cases}
\Sigma_{i=1}^{i=n} y_{2i}y_{2i-1} & p= 2, r=0 \\
\Sigma_{i=1}^{i=n}(y^{\iota p^r}_{2i-1}y_{2i}-
y_{2i-1} y^{\iota p^r}_{2i}) & {\rm otherwise },
\end{cases}
$$
with $\iota=1$ for $p=2$, and $\iota =p$ for $p>2$. Let $A_0$ be
an element of $\mathcal A$. We have the following theorem.

\begin{thm}[Green--Minh \cite{GM}] \quad
\label{thm:thm3.2}
\begin{enumerate}
\item
In $H^*(P)$ for $0 \leq r \leq n-1$,
$$
\kappa_{n,r}=Q_r(P/A_0)- \sum_{x \in {\P}
H^1(P/A_0)}\chi_{r,x}.
$$
\item
There exist
$f_{n,0}, \ldots , f_{n,n-1} \in  H^*(E)$, viewed as elements of
$H^*(P )$ via the inflation map, such that
$$
z_n^{(n)}+\sum_{i=0}^{n-1} (-1)^{n-i}z_n^{(i)}f_{n,i}=0,
$$
and, for every $A \in \mathcal A$,
$\res_A^P(\kappa^\iota_{n,r})=\res^P_A(f_{n,r})$, 
$0 \leq r \leq n-1$.
\item There exist $h_i$, $0 \leq  i \leq n-1$, and a unique 
$\eta$ of $H^*(E)$ such that
$$
z_n^{(n-1)}=y_{2n}\eta+\sum_{i=0}^{n-2} h_iz_n^{(i)},
$$
and in $h^*(P)$, $\chi_{n-1,x_{2n}}^{\iota}=-\inf(\eta^{p-1})$.
Furthermore, for all $0 \leq r \leq n-1$ and all 
$\phi \in {\P} H^1(E)$, 
$\chi^\iota_{r,\phi} \in \im(\inf^E_P)$, as elements of
$h^*(P)$.
\end{enumerate}
\end{thm}

By Quillen \cite{Q} it is known that, for $p = 2$, all the
$\kappa_{n,r}$ and $\chi_{r,\phi}$ belong to  $\mathcal T$.
For $p>2$, it follows that the above theorem that all 
$p^{\rm th}$--powers of the 
$\kappa_{n,r}$ and $\chi_{r,\phi}$ belong
to $\mathcal T$.
In fact, by setting
$$  
\varphi=   
\begin{cases}
y^p_{2n-1}-y_{2n-1}y_{2n}^{p-1}  & p> 2 \\
y_{2n-1} & p=2,
\end{cases}
$$
we have the following corollary.
\begin{cor} 
\label{cor:cor3.3}
In $h^*(P)$,
\begin{enumerate}
\item $\kappa_{n,r}^\iota=f_{n,r}$, $0 \leq r \leq n-1$
\item
$$
\eta=(-1)^{n-1}\left[\varphi f_{n-1,0}+\sum_{i=1}^{n-1}
(-1)^i\left(y^{\iota p^i}_{2n-1}-y_{2n-1}y_{2n}^
{\iota p^i-1}\right)f_{n-1,i}\right]
$$
\item with $K=\ker (x_{2n})$,
$$  
\eta=  
\begin{cases}
\tr^K_P(\zeta_{n-1}) & p= 2  \\
(-1)^{n-1}\tr^K_P \left[\zeta^{p-1}_{n-1}\sum_{i=0}^{n-1} (-1)^i
y^{p^i}_{2n-1}\kappa_{n-1,i}\right] & p>2;
\end{cases}
$$
\item for $0 \leq r \leq n-1$,
$$
\chi^\iota_{r,x_{2n}}=-f_{n-1,r}\left[\varphi
f_{n-1,0}+\sum_{i=1}^{n-1}(-1)^i\left(y^{\iota
p^i}_{2n-1}-y_{2n-1}y^{\iota p^i-1}_{2n}\right)f_{n-1,i}
\right]^{p-1} .
$$
\end{enumerate}
\end{cor}

\begin{proof}
Part (i) follows from \fullref{thm:thm3.2} (ii), by noting
that the restriction map from $h^*(P)$ to 
$\prod_{A \in \A} H^*A$ is injective.

We have, by \fullref{thm:thm3.2},
\begin{align*}
z^{(n-1)}_n & = z^{(n-1)}_{n-1}
+\left(y_{2n-1}^{\iota p^{n-1}}y_{2n}-
y_{2n-1}y_{2n}^{\iota p^{n-1}}\right) \\
& = \left(y_{2n-1}^{\iota p^{n-1}}y_{2n}-
y_{2n-1}y_{2n}^{\iota p^{n-1}}\right)
+ (-1)^n\sum_{i=0}^{n-2}(-1)^i z^{(i)}_{n-1}f_{n-1,i} \\
& = \left(y_{2n-1}^{\iota p^{n-1}}y_{2n}-
y_{2n-1}y_{2n}^{\iota p^{n-1}}\right) \\
&\qquad+ (-1)^n\biggl[\sum_{i=0}^{n-2}(-1)^i z^{(i)}_{n}f_{n-1,i}
-y_{2n}\varphi f_{n-1,0} \\[-2ex]
&\hspace{100pt}-\sum_{i=1}^{n-2}(-1)^i\left(y_{2n-1}^{\iota p^{i}}y_{2n}
-y_{2n-1}y_{2n}^{\iota p^{i}}\right)f_{n-1,i}\biggr] \\
& = (-1)^n \left[
\sum_{i=0}^{n-2}(-1)^iz^{(i)}_{n}f_{n-1,i}-y_{2n}X
\right],
\end{align*}
with 
$$
X=\varphi f_{n-1,0}+\sum_{i=1}^{n-1}(-1)^i
(y_{2n-1}^{\iota p^{i}}-y_{2n-1}y_{2n}^{\iota p^{i}-1})f_{n-1,i}. 
$$
So $\eta = (-1)^{n+1}X$; (ii) and (iv) are proved.

Pick an element $A \in \A$. By \cite[Lemma 7.1]{GM} and its 
proof, we have
$$
\res^P_A\eta = 
\begin{cases}
-V(B',y_{2n-1})^\iota  & A \subseteq K\\
0  & \rm{otherwise.} 
\end{cases}
$$
Set
$$  
Y = 
\begin{cases}
\tr^K_P(\zeta_{n-1})& p= 2 \\
(-1)^{n-1}\tr^K_P[\zeta^{p-1}_{n-1}
\sum_{i=0}^{n-1}y^{p^i}_{2n-1} \kappa_{n-1,i}] & p>2.
\end{cases}
$$
If $A \not \subset K$, then $\res^P_A(Y)=0$, by the Mackey
formula. Suppose $A \subset K$. Set $B=A \cap P_{n-1}$. For 
$p>2$,
we have
\begin{eqnarray*}
\res^P_A(Y) & = & (-1)^{n-1} \sum_{g \in P/K} {}^g
\sum_{i=0}^{n-1}\res^K_A 
((-1)^{i}y^{p^i}_{2n-1}\ka_{n-1,i}\zeta^{p-1}_{n-1}) \\
& = & (-1)^{n-1}\sum_{g \in P/K }
{}^g\sum_{i=0}^{n-1}(-1)^{i}y^{p^i}_{2n-1} 
Q_i(B')\zeta^{p-1}_{n-1} \\
& = & (-1)^{n-1}\sum_{i=0}^{n-1} (-1)^{i} 
y^{p^i}_{2n-1}Q_i(B')\sum_{g \in P/K} {}^g\zeta^{p-1}_{n-1}
\end{eqnarray*}
since the $y_{2n-1}$ and the $Q_i(B')$ are invariant under the
action of $P/K$. Thus
$$   
\res^P_A(Y)=   
\begin{cases}
\res^P_A (\tr^K_P(\zeta_{n-1})) & p= 2  \\
V(B',y_{2n-1}) \res^P_K (\tr^K_P(\zeta_{n-1}^{p-1})) & p>2.
\end{cases}
$$
Following \cite[Proposition 4.4]{GM} we have
$$
\res^P_A (\tr^K_P(\zeta_{n-1}^{p-1}))=-V(B',y_{2n-1})^{p-1}.$$ So
$\res^P_A(Y)=-V(B',y_{2n-1})^\iota$.

Since $\eta - Y$ restricts trivially to every element of $\A$, it
follows that $\eta \doteqdot Y$. 
\end{proof}

\begin{prop}
\label{prop:prop3.4}
For $0 \leq r \leq n-1$, 
$$ 
\ka_{n,r} \doteqdot - \sum_{x \in \P H^1(P) } \chi_{r,x}. 
$$
\end{prop}

\begin{proof}
Let $A$ be an element of $\A$. There exist exactly
$\smash{\frac{p^n-1}{p-1}}$ elements of $\P H^1(P)$ of which the kernel
contains $A$. The subset of those elements is nothing but 
$\P H^1(P/A)$.

Let $x$ be an element of $\P H^1(P)$. It is clear that
$\res^P_A(\chi_{r,x})=0$ if $x \not \in \P H^1(P/A)$.

Hence
$$
\res^P_A\Bigl(\sum_{x \in \P H^1(P)}  \chi_{r,x}\Bigr)=
\res^P_A\Bigl(\sum_{x \in \P H^1(P/A)} \chi_{r,x}\Bigr).
$$
Therefore, by \cite [Theorem 5.2]{GM} 
\begin{eqnarray*}
\res^P_A\Bigl(\sum_{x \in \P H^1(P)} \chi_{r,x}\Bigr)
& = & \res^P_A\Bigl(Q_r(P/A) - \ka_{n,r}\Bigr) \\
& = & -Q_r(A'),
\end{eqnarray*}
since $\res^P_A(Q_r(P/A))=0$ and $\res^P_A(\ka_{n,r})=Q_r(A')$.
The proposition follows.
\end{proof}

We are now going to obtain characteristic classes of $P$ by using
the Evens norm map. We first need the following lemma.
\begin{lem}
\label{lem:lem3.5}
Fix a generator $e$ of $Z$. Let $H= A \cap B$ with
$A,B \in \mathcal A$ and let $(h_1, \ldots , h_ k, e)$ be a basis
of $H$. Then there exist elements $g_1, \ldots , g_k$ of $P$
satisfying
\begin{enumerate}
\item $[g,g_i] = 1$, 
$
[g_i,h_j] =  
\begin{cases}
1 & i \neq j  \\
e & i=j
\end{cases}
$ and $1 \leq i,j \leq k$.
\item
$P= \coprod_{g \in G} AgB$ is a double coset decomposition of
$P$ with $G=\langle g_1,\ldots,g_k \rangle$. 
\end{enumerate}
\end{lem}

\begin{proof}
The existence of the $g_j$ satisfying (i) follows
from \cite{M1}. Assume that $agb = a' g' b'$ with $a, a' \in A$,
$b,b' \in B$, $g,g' \in G$. It follows that 
$[g,h_i]= [g', h_i]$ and $1\leq i\leq k$, hence $g=g'$. 
As $|G|=p^k$, (ii) is obtained. 
\end{proof}

The following notation will be used. Let $C$ be the cyclic group
of order $p$ and fix a generator $u$ of $H^1(C)$ (resp.\ $H^2(C)$)
for $p = 2$ (resp.\ $p > 2$). Set $\Gamma = P \times C$. If $H$ is
a subgroup of $P$, every element of $H^*(H)$ 
(resp.\ $H^*(C)$) can
be considered as an element of $H^*(H \times C)$. We have the
following lemma.

\begin{lem}
\label{lem:lem3.6}
Let $A$, $B$ be elements of $\mathcal A$ and let 
$v \in  H^1(A \times C)$ (resp.\ $H^2(A \times C)$) for $p=2$ 
(resp.\ $p>2$). Assume that 
$\res^{A \times C}_{Z \times C}(v)=\mu \gamma+\lambda u$ with 
$\mu, \lambda \in \F_p$, then
$$
\res^{\Gamma}_{B \times C} 
{\mathcal N}_{A \times C \rightarrow \Gamma}(v)= 
\mu V(B',\gamma)+\lambda V(B',u).
$$
\end{lem}

\begin{proof}
Let $P= \cup_{g \in G} AgB$ be the double coset
decomposition of $P$ given in \fullref{lem:lem3.5}. 
Set $H = B \cap A$. We have \\
\begin{align*}
\res^\Gamma_{B \times C} \N_{A \times C \sra \Gamma}(v)
& = \prod_{g\in G} \N_{(B \cap A^g) \times C \sra B \times C}
\res^{A^g \times C}_{(B \cap A^g) \times C}({}^gv) \\
& = \prod_{g \in G} \N_{H \times C \sra B \times C} \circ
\res_{H \times C}^{A \times C}({}^gv) & \hbox{($A$ is normal)} \\
& = \N_{H \times C \sra B \times C} \biggl( \prod_{g \in G} 
\res_{H\times C}^{A \times C}({}^gv)\biggr) 
& \hbox{(\fullref{cor:cor2.3})} \\
& =  \N_{H \times C \ra B \times C} \biggl( \prod_{g \in
G}(\mu^g\gamma+\lambda u)\biggr) \\
& = \N_{H \times C \ra B \times C} \bigl( V(H',\mu \gamma +\lambda 
u)\bigr) & \hbox{(\fullref{lem:lem3.5})}\\
& = \N_{Z \times C \ra B \times C} ( \mu \chi + \lambda u)
& \hbox{(\fullref{cor:cor2.2})}\\
& =  \mu V(B',\gamma) + \lambda V(B',u) & \hbox{(\fullref{cor:cor2.3})}
\end{align*}
as required.
\end{proof}

The following shows that characteristic classes of $P$ can be
obtained by means of the Evens norm map.
\begin{thm}
\label{thm:thm3.7}
Let
$A$ be an element of $\A$ and let $v$ be an. element of $H^1(A)$
(resp.\ $\beta H^1(A)$) for $p = 2$ (resp.\ $p > 2$) satisfying
$\res^A_Z(v)=\gamma$.  Set
$$
\zeta_{A,v} = 
\N_{A \times C\ra \Gamma}(v + u) - 
\N_{A\times C\ra \Gamma}(v).
$$
As elements of $h^*(P )$ then
$$ 
\zeta_{ A,v} = (-1)^n\sum_{s=0}^{n}(-1)^s \kappa_{n,s}u^{p^s}.
$$
\end{thm}

\begin{proof}
For every $B \in  \A$, by \fullref{lem:lem3.6} we have
\begin{eqnarray*}
\res_{B\times Z}^\Gamma(\zeta_{A,v}) &=& 
\res_{B\times Z}^\Gamma
\N_{A\times Z \ra \Gamma}(v + u) - \res_{B\times Z}
\N_{A\times Z \ra \Gamma}(v )\\
& =& V (B' , \gamma + u) - V (B' ,\gamma ) \\
&=& V (B' , u),
\end{eqnarray*}
since $V (B' , X)$, as a function on $X$, is additive. By 
\fullref{thm:thm3.1} (ii),
$$
\res_{B\times Z}^\Gamma(\zeta_{A,v}) = \res_{B\times
Z}^\Gamma[(-1)^n \sum_{s=0}^{n}(-1)^s\kappa_{n,s}u^{p^s}].
$$
So $\zeta_{A,v}\doteqdot (-1)^n
\sum_{s=0}^{n}(-1)^s\kappa_{n,s}u^{p^s}$.
\end{proof}

\begin{rem} 
Write $\N_{A\times C\ra \Gamma} = \N$. It follows
from the above theorem and from \fullref{cor:cor2.2}
that
$$ 
\N(v + u) - \N(v) - \N(u)\doteqdot (-1)^n 
\sum_{s=0}^{n}(-1)^s[\kappa_{n,s} - Q_s(P/A)]u^{p^s}.
$$
According to \fullref{prop:prop2.1}(iv), the $\kappa_{n,s}$ can be
expressed as sums of transfers from maximal subgroups of $P
\times C$. Such formulae are the ones given in \fullref{thm:thm3.2}.
\end{rem}

Let $a_1, a_2, \ldots , a_{2 n-1} , a_{2 n}$ be elements of $P$
satisfying $x_i(a_j) = \delta_{ij}$ with $\delta_{ij}$ the
Kronecker symbol, $1 \leq i,j \leq 2n$. Suppose that $K$ is a
maximal subgroup of $P$ given by $K = \ker (x_{2n})$. So 
$K \cong P_{n-1}\times \langle a_{2n-1} \rangle 
\cong P_{n-1}\times \Z/p$.
Write $y_{2n}=y$, $\N_{K\times C\ra \Gamma} = \N$, and, for $0
\leq r \leq n-1$, $\chi_{r,x_{2n}}=\chi_r$. Define
$$
\theta_{n-1, r}=\res^P_K(\chi_r) \in H^*(K),
$$
and
\begin{align*}
\theta_{n-1}=u^{p^n}+\sum_{r=0}^{n-1}
\left[(-1)^{n-r}u^{p^r}(-\theta_{n-1,r}+\kappa^p_{n-1,r-1})\right]\in
H^*(K \times C)
\end{align*}
with the convention that $\kappa_{n-1,-1}=0$. 

\begin{thm}
\label{thm:thm3.8}
As elements of $h^*(P )$,
\begin{align*}
(-1)^n \N(\theta_{n-1})
& =  \sum_{s=0}^{n}(-1)^s\kappa^p_{n,s}u^{p^{s+1}}\\
(-1)^r \N(\kappa_{n-1,r}^\iota)
& =  \sum_{i=r}^{n-1}(-1)^i\kappa^{\iota
p}_{n-1,i}y^{\iota(p^{i+1}-p^{r+1})} \\
& \qua - \biggl(\sum_{i=r+1}^{n-1}(-1)^i\kappa^{\iota
}_{n-1,i}y^{\iota(p^{i}-p^{r+1})}\biggr)
\biggl[\sum_{i=0}^{n-1}(-1)^i\kappa^{\iota}_{n-1,i}
y^{\iota p^{i}}\biggr]^{p-1},
\end{align*}
for $0 \leq r \leq n-2$.
\end{thm}

\begin{proof}
For convenience, write $\kappa_{n-1 ,r}=\kappa_r$ for 
$0 \leq r \leq n-1$. Let $A$ be an element of $\A$ and set 
$X = \N(\theta_{n-1})$ and $Y_r = \N (\kappa^\iota_{n-1,r})$.
Let
$$Z_r = \sum_{i=r}^{n-1}(-1)^i\ka^{\iota p}_{i}
y^{\iota(p^{i+1}-p^{r+1})}
- \biggl(\sum_{i=r+1}^{n-1}\!\!(-1)^i
  \ka^{\iota}_{i}y^{\iota(p^{i}-p^{r+1})}\!\biggr)
\biggl[\sum_{i=0}^{n-1}(-1)^i\ka^{\iota}_{i}y^{\iota p^{ i}}
\biggr]^{p-1}$$
for $0 \leq r \leq n-2$.  Consider the following cases:

\textbf{Case 1} $A \subset K$\qua By
setting $B = A  \cap P_{n-1}$, we have $A = B \times \langle 
a_{2n-1}\rangle$.
So $$ \res^\Gamma_{K\times C}(X) = \prod_{x \in
\langle a_{2n}\rangle}{}^x\theta_{n-1} .
$$
As the $\theta_{n-1}$ belong to $\im(\res^P_K)$, they are
invariant under the action of $a_{2n}$. Hence
\begin{align*}
\res_{A\times C}^\Gamma(X)& = \prod_{x \in
\langle a_{2n}\rangle }\res_{A\times C}^{K \times C} 
({}^x\theta_{n-1}) \\
&=\res_{A\times C}^{K \times C} (\theta_{n-1}^p) \\
&=u^{p^{n+1}} {+} 
\biggl[ 
\sum_{r=0}^{n-1}(-1)^{n-r}u^{p^r}\biggl(Q_r(B' )V
(B' , y_{2 n-1})^{p-1}
{+} Q^p_{r-1}(B')\biggr)
\biggr]^p
\\
&=V (A' , u)^p \\
& =\res_{A\times C}^\Gamma \biggl((-1)^n 
\sum_{s=0}^{n}(-1)^s\ka^p_{n,s}u^{p^{s+1}}\biggr).
\end{align*} 
Also, for $0 \leq r \leq n - 2$,
\begin{multline*}
$$\res_{A\times C}^\Gamma(Y_r)
= \prod_{x \in \langle a_{2n}\rangle }
\res_{A\times C}^{K \times C} 
\left({}^x\ka^\iota_{n-1,r}\right) \\
= \left[\res_{A\times C}^{K \times C} 
(\ka_{n-1,r})^\iota \right]^p
= Q^{\iota p}_r(B')
=  (-1)^r \res^\Gamma_{A \times C}(Z_r).
\end{multline*}
\textbf{Case 2} $A \not \subset K$\qua By setting $H = K  \cap A$, we have
\begin{align*}
\res_{A\times C}^\Gamma(X)
& =  
\N_{H \times C \sra A \times C} \res_{H\times C}^{K \times C} (X) \\
&= \N_{H \times C \sra A \times C} (V(H',u)^p)
\\
&= V(A',u)^p
\\
&= 
\res_{A\times C}^\Gamma \biggl((-1)^n\sum_{s=0}^{n}(-1)^s 
\ka^p_{n,s}u^{p^{s+1}}\biggr).
\end{align*}
and
\begin{align*}
\res_{A\times C}^\Gamma(Y_r)& = \N_{H \times C \sra A \times 
C} \res_{H\times C}^{K \times C} (Y_r) \\
&= \N_{H \times C \sra A \times C} (Q^\iota_r(H/Z)) \\
&= (-1)^r\res_{A\times C}^\Gamma(Z_r). 
\end{align*}
This completes the proof.
\end{proof}

Formulae relating the $\ka^\iota_{n,r}$ to such classes of fewer
variables are given by the following corollary.

\begin{cor}
\label{cor:cor3.9}
For $0 \leq r \leq n-1$, as elements of $ h^*(P)$,
\begin{multline*}
\ka^\iota_{n,r}=\ka^{\iota
p}_{n-1,r-1}+\ka^\iota_{n-1,r}\biggl[\sum^{n-1}_{i=0}  (-1)^i
\ka^\iota_{n-1,i}y^{\iota p^i}\biggr]^{p-1}\\[-1ex]
+\biggl[\ka^\iota_{n-1,0}\varphi+\sum_{i=1}^{n-1}
(-1)^i\kappa^\iota_{n-1,i}(y^{\iota p^i}_{2n-1}-y_{2n-1}y^{\iota
p^i-1}_{})\biggr]^{p-1}.
\end{multline*}
\end{cor}

\begin{proof}
By \fullref{cor:cor2.3} we have
\begin{align*}
\N(\theta_{n-1}) & = \N(u^{p^n})+\sum_{r=0}^{n-1}
\N\left[\left(-1)^{n-r}u^{p^r}
(-\theta_{n-1,r}+\ka^p_{n-1,r-1}\right)\right]\\[-1ex]
&= V(y,u)^{p^n}+\sum_{r=0}^{n-1} 
(-1)^{n-r}V(y,u)^{p^r}\left[\N(-\theta_{n-1,r})
+\N(\ka^p_{n-1,r-1})\right]  \\
&= (u^p -uy^{p-1})^{p^n}  \\[-1ex]
& \qquad + \sum_{r=0}^{n-1} (-1)^{n-r}(u^p 
-uy^{p-1})^{p^r}\left[\N(-\theta_{n-1,r})
+\N(\ka^p_{n-1,r-1})\right]\\[-1ex]
&= u^{p^{n+1}}+\sum_{r=0}^{n-1} 
(-1)^{n-r}u^{p^{r+1}}\biggl[\N(-\theta_{n-1,r})
+\N(\ka^p_{n-1,r-1}) \\[-2ex]
& \hspace{120pt} +y^{(p-1)p^{r+1}}
\left(\N(-\theta_{n-1,r+1}) +\N(\ka^p_{n-1,r})\right)\biggr] \\[-1ex]
& \qua +(-1)^nuy^{p-1}\N(-\theta_{n-1,0}).
\end{align*}
By the Frobenius formula, the cup-product of $\chi_r$ with each 
of
$x_{2n}$ , $y_{2n}$ vanishes. As the transfer commutes with
Steenrod operations, we have, by \fullref{prop:prop2.4} and
\fullref{thm:thm3.8},
\begin{align*}
&u^{p^{n+1}}+ \sum_{r=0}^{n-1}(-1)^{n-r}
\ka^p_{n,r}u^{p^{r+1}} =
\N(\theta_{n-1})  \\[-1ex]
& = u^{p^{n+1}}+ \sum_{r=0}^{n-1} (-1)^{n-r}
u^{p^{r+1}}\biggl[-\chi^p_{r} +\N(\ka^p_{n-1,r-1}) \\[-2ex]
& \hspace{200pt} +y^{(p-1)p^{r+1}}
\left(-\chi^p_{r+1}+\N(\ka^p_{n-1,r})\right)\biggr] \\[-1ex]
& = u^{p^{n+1}}+
\sum_{r=0}^{n-1}
(-1)^{n-r}u^{p^{r+1}}\left[-\chi^p_{r}+\N(\ka^p_{n-1,r-1}) 
+y^{(p-1)p^{r+1}}\N(\ka^p_{n-1,r})\right]. 
\end{align*}
Therefore
$$\sum_{r=0}^{n-1} (-1)^{n-r} \ka^{\iota}_{n,r}
  u^{\iota p^r}\!=\sum_{r=0}^{n-1} (-1)^{n-r}u^{\iota p^r}
  \biggl[-\chi^\iota_{r}{+}\N(\ka^\iota_{n-1,r-1})
  {+}y^{\iota (p-1)p^{r}}\!\N(\ka^\iota_{n-1,r})\biggr].$$
Hence
$$
\ka^{\iota}_{n,r}=-\chi^\iota_{r}+\N(\ka^\iota_{n-1,r-1})+
y^{(p-1)\iota p^{r}} \N(\ka^\iota_{n-1,r}).
$$
Since
$$
\N(\ka^\iota_{n-1,r-1})+y^{(p-1)\iota p^{r}} \N(\ka^\iota_{n-1,r}) 
\doteqdot \ka^{\iota
p}_{n-1,r-1}+\ka^{\iota }_{n-1,r}\left[\sum_{i=0}^{n-1} (-1)^i
\ka^{\iota }_{n-1,i}y^{\iota p^{i}}\right]^{p-1},
$$
by \fullref{thm:thm3.8}, we obtain
$$
\kappa^\iota_{n,r} \doteqdot - \chi^\iota_r+ \kappa^{\iota
p}_{n-1,r-1}+\kappa^{\iota }_{n-1,r}\left[\sum_{i=0}^{n-1}
(-1)^i\kappa^\iota_{n-1,i} y^{\iota p^i}\right]^{p-1} .
 $$
The corollary follows from \fullref{cor:cor3.3}.
\end{proof}

\section[The subring]{The subring 
$\F_p[\ka_{n,0}, \ldots ,\ka_{n,n-1}] \cap \T$}
\label{sec:sec4}

In this section, $p$ is supposed to be an odd prime. It was 
proved by Green and Leary \cite{G,GL} that $\ka^s_{n,0} \in \T$,
provided that $s \geq 2^n$, or $s \geq 2$ and $n\leq 2$. This result
can be sharpened as follows. Let $ \Re_n$
be the set consisting of sequences $R=(r_0,r_1,\ldots,r_{n-1})$
of non-negative integers. For
$R= (r_0 , \ldots  , r_{n-1}) \in {\Re}_n$ and for $m>0$, set
\begin{align*}
s_R & =\sum_{i \geq 0} r_i, \\
\ka^R_m& = \begin{cases}
\Pi_{i=0}^{i=m-1}\ka^{r_i}_{m,i}  & m \leq n    \\
\Pi_{i=0}^{i=n-1}\ka^{r_i}_{m,i} & m > n.
\end{cases}
\end{align*}
The main purpose of this section is to prove the following theorem.
\begin{thm} 
\label{thm:thm4.1}
Let $R = (r_0, \ldots , r_{n-1})$ be an element of 
${\Re}_n$. As
an element of $h^*(P )$, $\ka^R_n$ belongs to $\T$ if and only if
one of the following conditions is satisfied:\nl
\begin{tabular}{ll}
$(R_1)$ & $r_0 \geq 2$; \\
$(R_2)$ & $r_0=0$ and all the $r_i$ terms with $i >0$, are multiples of $p$.
\end{tabular}
\end{thm}

The rest of the section is devoted to the proof of the theorem. 

\begin{proof}
By \fullref{cor:cor3.3}, 
$\ka^R_n \in \T$ if $R$ satisfies $(R_2)$. We
shall prove the following proposition.
\begin{prop}
\label{prop:prop4.2}
If $R \in {\Re }_n$ satisfies $(R_1)$, then $\ka^R_n \in \T$.
\end{prop}
By \cite{GL,L}, the proposition holds for $n = 1$.
Suppose inductively that it holds for $n - 1$. Set 
$K = \ker (x_{2n}) = P_{n-1} \times \Z /p$ and 
$\T '= \im (\inf_K^{K/Z}) + \sqrt 0$. 
Write $w = \res^P_K(y_{2n-1})$ and 
$\psi=(-1)^{n-1} \sum_{j=0}^{n-1} (-1)^j \ka_{n-1,j}w^{p^j}$. We have
$$
\res^P_K(\ka_{n,j}) \doteqdot
\ka_{n-1,j-1}^p+\ka_{n-1,j}\psi^{p-1} , \qua 0 \leq j \leq n.
$$
So, for every element $R \in \R_n$, as elements of $h^*(K)$,
\begin{align}
\res_K^P(\ka^R_n) &= \prod_{j=0}^{n-1}
\left[\ka^p_{n-1,j-1}+\ka_{n-1,j}\psi^{p-1}\right]^{r_j} \nonumber \\
&=
\ka^{r_0}_{n-1,0}\psi^{r_0(p-1)}\prod_{j=1}^{n-1}\left[\ka^p_{n-1,j-1}
+\ka_{n-1,j}\psi^{p-1}\right]^{r_j} \nonumber \\
&=\ka^{R}_{n-1}\psi^{(p-1)s_R} + \sum_{r_0 \leq t  < s_R} 
\rho_t
\psi^{(p-1)t} 
\label{eqn:eqn2}
\end{align}
with $\rho_t \in h^*(P_{n-1})$.

\begin{lem}
\label{lem:lem4.3}
Let $S = (s_0, \ldots, s_{n-1})$ be an 
element of $\R_n$ with $s_0\geq 1$, and let $x$
be a non-zero element of $H^1(P )$. Then
$$
\ka^S_n \chi_{0 ,x} \in  \T .
$$
\end{lem}

\begin{proof}
Without loss of generality, we may assume that $x =x_{2n}$. 
So $K = \ker (x)$. Since $s_0 \geq 1$, by \eqref{eqn:eqn2}, we have
\begin{equation}
\res^P_K (\ka^S_n) \doteqdot \sum_{U \in \U} \ka^U_{n-1}w^{t_U}
\psi
\label{eqn:eqn3}
\end{equation}
with $\U$ a subset of
$$
\{R = (r_0,\ldots,r_{n-2}) \in \R_{n-1}|r_0 \geq  1\}.
$$
Let $U = (u_0,\ldots, u_{n-2})$ be an element of $\U$. Since
$$
\ka^U_{n-1}\ka_{n-1,0}= \ka_{n-1,0}^{u_0+1} \prod_{i=1}^{n-2}
\ka^{u_i}_{n-1,i}
$$
and $u_0+1 \geq 2$, it follows from the inductive hypothesis that
$\ka^U_{n-1}\ka_{n-1,0}$, and hence
$\ka^U_{n-1}\ka_{n-1,0}w^{t_U}$ belong to $\T'$. So, via the
inflation map, $\ka^U_{n-1}\ka_{n-1,0}w^{t_U}$ belongs to $\T$.

We then have as elements of $h^*(P)$,
\begin{align*}
\ka^S_n\chi_{0,x} &=\ka^S_n
\tr^K_P(\ka_{n-1,0}\zeta^{p-1}_{n-1})
\\
&=\tr^K_P(\res^P_K(\ka^S_K).\ka_{n-1,0}\zeta^{p-1}_{n-1}) && \text{by Frobenius formula} \\
&=\sum_{U \in \U} \tr^K_P(\ka^U_{n-1}\ka_{n-1,0}w^{t_U}
\psi\zeta^{p-1}_{n-1}) && \mathrm{by~\eqref{eqn:eqn3}}\\
&=\sum_{U \in \U}\ka^U_{n-1}\ka_{n-1,0}w^{t_U}.\tr^K_P
(\psi\zeta^{p-1}_{n-1}),
\end{align*}
which implies  $\ka^S_n\chi_{0,x} \in \T$, by \fullref{cor:cor3.3} (iii).
\end{proof}

\begin{proof}[Proof of \fullref{prop:prop4.2}]
Let $R =(r_0, \ldots,r_{n-1})$ be
an element of $\R_n$. By \fullref{cor:cor3.3}(i), $\ka^R_n \in \T$ if
$R$ satisfies $(R_2)$. Suppose that $r_0 \geq 2$. Set
$S = (r_0 - 1, r_1, \ldots , r_{n-1})$. We then have
$$
\ka^R_n=\ka^S_n\ka_{n,0}=-\sum_{x \in \P H^1(P)}
\ka^S_n\chi_{0,x},
$$
by \fullref{prop:prop3.4}. Since $r_0-1 \geq 1$ by \fullref{lem:lem4.4}
$\ka^S_n\chi_{0,x} \in \T$, for every $x \in \P H^1(P)$;
so $\ka^R_n \in \T$. The proposition is proved.
\end{proof}

Consider $\psi$, and also the right hand side of \eqref{eqn:eqn2}, as
polynomials with variable $w$ and with coefficients in
$h^*(P_{n-1})$. We have the following lemma.
\begin{lem}
\label{lem:lem4.4}
Let $R =(r_0, \ldots,r_{n-1})$ be an element of $\R_n$ with $s_R
\not = 0$ mod $p$. Then for $0 \leq i \leq n-2$,
\begin{enumerate}
\item
$\res^P_K(\ka^R_n) \doteqdot
s_R(-1)^{i+n}\ka^R_{n-1}\ka_{n-1,i}w^{p^{n-1}[(p-1)s_R-1]+p^i}+$ other terms;
\item
$\ka^R_{n-1}\ka_{n-1,i} \in \T$ if $\ka^R_n \in \T$.
\end{enumerate}
\end{lem}

\begin{proof}
For $t < s_R$, ${\rm deg}(\psi^{(p-1)t}) \leq
p^{n-1}(p-1)s_R-p^n +p^{n-1}$; hence
$$
{\rm deg}(\psi^{(p-1)t}) < {\rm min}(p^{n-1}(p - 1)s_R - 1,
p^{n-1}[(p - 1)s_R - 2] + p^i).
$$
So (i) follows from \eqref{eqn:eqn2} and the fact that
\begin{align*}
\psi^{(p-1)s_R} & = [ \sum_{i=0}^{n-1}(
-1)^{i}\ka_{n-1,i}w^{p^i}]^{(p-1)s_R} \\
& =
-s_R\sum_{i=0}^{n-2}(-1)^{i+n-1}
\ka_{n-1,i}w^{p^{n-1}[(p-1)s_R-1]+p^i} + \text{other terms.}
\end{align*}

Write
$$
\res^P_K (\ka^R_n) = \sum_{i \geq 0} \rho_i  w^i
$$
with $\rho_i \in h^*(P_{n-1})$. If $\ka^R_n \in \T$, then
$\res^P_K (\ka^R_n)$ belongs to $\T'$, so all the $\rho_i$ lie
in $\T'$; (ii) is then a direct consequence of (i).
\end{proof}

The proof of the theorem is completed by \fullref{prop:prop4.2}
and the following.
\begin{lem} 
\label{lem:lem4.5}
If $\ka^R_n \in \T$   with $R = (r_0, \ldots , 
r_{n-1}) \in \R_n$, then
$R$ satisfies $(R_1)$ or $(R_2)$. 
\end{lem}

\begin{proof}
By Leary \cite{L}, the lemma holds for $n = 1$. Assume
that it holds for $n - 1$.

Suppose that $\ka^R_n \in \T$ with  $R = (r_0, \ldots , r_{n-1})$
and $r_0 < 2$. It follows that 
$\zeta = \res^P_K( \ka^R_n) \in \T'$.
Consider $\zeta$ as a polynomial with variable $w$ and with
coefficients in $h^*(P_{n-1})$. By \eqref{eqn:eqn2}, we have 
\begin{align*}
\zeta &=\ka^R_{n-1} \psi^{(p-1)s_R}+\sum_{0 \leq t <
s_R}\rho_t\psi^{(p-1)t} \\
&=\ka^R_{n-1}w^{p^{n-1}(p-1)s_R}+ \text{other terms}, & \text{by \eqref{eqn:eqn3}}
\end{align*}
 which implies $\ka^R_{n-1} \in \T'$. By the induction 
hypothesis, $r_0 = 0$ and
$r_1,\ldots ,r_{n-2}$ are multiples of $p$. So $s_R=r_{n-1}$ mod
$p$. If $s_R \not = 0$ mod $p$, it follows from 
\fullref{lem:lem4.5} that
$$
\ka^{r_1}_{n-1,1}\ldots\ka^{r_{n-3}}_{n-1,n-3}\ka^{r_{n-2}+1}_{n-
1,n-2}=\ka^{}_{n-1,n-2}\ka^{R}_{n-1}
\in \T'
$$
which contradicts the induction hypothesis, since $r_{n-2} = 0$
mod $p$ implies $r_{n-2}+1 \not = 0$ mod $p$. So $s_R = 0$ mod
$p$, hence $r_{n-1} = 0$ mod $p$. The lemma follows.
\end{proof}

This completes the proof of \fullref{thm:thm4.1}.
\end{proof}

Let $x$ be a non-zero element of $H^1(P )$. By \fullref{thm:thm3.2}(iii),
there exists a unique $\eta_x \in H^*(E)$ such that, as elements
of $H^*(E)/(z_n^{(1)} , \ldots , z_n^{(n-2)})$,
\begin{equation}
z_n^{ (n-1)} = \begin{cases} \eta_x \beta(x) & 
\text{$p$ odd,} \\
\eta_x x & \text{$p=2$.} \end{cases}
\label{eqn:eqn4}
\end{equation}
Note that $H_x = \ker (x)$ can be identified with $P_{n-1} \times
\Z/p$. Pick a non-zero element $u_x$ of $H^1(P )$ satisfying $$ 0
\not = \res^P_{H_x} (u_x) \in \ker (\res^{H_x}_{P_{n-1}}) \qua \text{set } 
v_x = \begin{cases}
u_x & p=2 \\ \beta(u_x) & \text{$p$ odd}, 
\end{cases}
$$
and define $\psi_x = (-1)^{n-1} \sum_{j=0}^{n-1}  (-1)^j
\ka_{n-1,j} v_x^{p^j}$.

 Let $R$ be an element of $\R_n$ satisfying $(R_1)$ or
$(R_2)$. By \fullref{thm:thm4.1}, $\ka^R_n$ and 
$\ka^R_{n-1}$ both belong
to $\T$. It is then interesting to find out a formulae relating
$\ka^R_n$ and $\ka^R_{n-1}$.
 If $R$ satisfies $(R_2)$, the formula can be derived from 
\fullref{cor:cor3.9}. In the case where
$R$ satisfies $(R_1)$ the formula follows from the next
corollary.

\begin{cor} 
\label{cor:cor4.6}
Let $R =(r_0, \ldots,r_{n-1})$ be an element of 
$\R_n$ with $r_0 \geq 2$.
Then, as elements of $\T$,
$$
\ka^R_n = -\ka^{r_0}_{n-1} \sum_{x \in \P H^1(P)} \eta_x
\psi_x^{(r_0-1)(p-1)-1} \prod^{n-1}_{j=1}[\ka^p_{n-1 ,j-1} +
\ka_{n-1,j} \psi_x^{p-1}]^{r_j}.
$$
\end{cor}

\begin{proof}
Set $S = (r_0 - 1, r_2, \dots , r_{n-1})$ and 
$U = (r_0 - 2, r_1, \ldots , r_{n-1})$. It follows from the proof of
\fullref{prop:prop4.2} that
\begin{align*}
\ka^R_n & = -\!\!\!\!\sum_{x\in  \P H^1(P)} \!\!\!\!\ka^S_n \chi_{0,x} \\
& = -\!\!\!\!\sum_{x \in \P H^1(P)} \!\!\!\!\ka^S_n \tr^{H_x}_P(\ka_{n-1,0}
\zeta_{n-1}^{p-1}) \\
&= -\!\!\!\!\sum_{x \in \P H^1(P)} \!\!\!\!\tr^{H_x}_P\left[\res^P_{H_x}
(\ka^S_n)\ka_{n-1,0}\zeta_{n-1}^{p-1}\right] \\
&= -\!\!\!\!\sum_{x \in \P H^1(P)} \!\!\!\!\tr^{H_x}_P\left[\res^P_{H_x}
(\ka^U_n)\ka^2_{n-1,0}\zeta_{n-1}^{p-1}\right] \\
&= -\!\!\!\!\sum_{x \in \P H^1(P)}\!\!\!\!\tr^{H_x}_P\biggl(\psi
\zeta_{n-1}^{p-1}\ka^{r_0}_{n-1,0}\psi_x^{(r_0-1)(p-1)-1}
%& \hspace{150pt}
\prod_{j=1}^{n-1} [\ka^p_{n-1 ,j-1}{+}\ka_{n-1,j}\psi^{p-1}_x]^{r_j}\biggr).
\end{align*}
Since $r_0 \geq 2$, it follows from \fullref{thm:thm4.1} that
$$
\rho_x = \ka^{r_0}_{n-1,0}\psi_x^{(r_0-1)(p-1)-1}
\prod_{j=1}^{n-1} \left[\ka^p_{n-1
,j-1}+\ka_{n-1,j}\psi^{p-1}_x\right]^{r_j}
$$
belongs to $\T$ , for any $x \in \P H^1(P)$. Hence
\begin{align*}
\ka^R_n &= - \sum_{x \in \P H^1(P)}  
\rho_x\tr^{H_x}_P(\psi_x
\zeta_{n-1}^{p-1}) \\
&= -\sum_{x \in \P H^1(P)} \rho_x \eta_x &\text{by \eqref{eqn:eqn4}}.
\end{align*}
This completes the proof.
\end{proof}

\section{Symplectic invariants}
\label{sec:sec5}

Recall that the symplectic group $ Sp_{2 n} = Sp_{2n} (\F_p)$ is
the group consisting of $E$  which preserve the nondegenerate
symplectic form $x_1x_2 + \ldots + x_{2n-1}x_{2n}$ of $H^2(E)$. 
Clearly $z_n^{(0)},\ldots,z_n^{(n-1)}$ belong to the
subring of invariants of  $ Sp_{2 n}$ in $\F_p[y_1, \ldots,
y_{2n}]$. According to a result of Quillen \cite{Q} for 
$p = 2$, and of Tezuka--Yagita \cite{TY} for $p > 2$, 
$$
\T =
\F_p[y_1,\ldots,y_{2n} ]/(z_n^{(0)} , \ldots , z_n^{(n-1)}).
$$
There is then an induced action of $ Sp_{2 n}$ on $\T$. Set
$$
\R'=\{R \in \R |R \text{ satisfies } (R_1) \text{ or } 
(R_2)\},
$$
and let $\R''$ be the subset of $\R_n$ consisting of elements $R 
=
(r_0, r_1,\ldots , r_{n-1})$ of $\R_n$ satisfying the following
two conditions:
\begin{itemize}
\item $0 \leq r_i \leq p-1$ for $i>0$
\item $r_0=3$, or $r_0=2$ and $r_1,\ldots,r_{n-1}$  are not 
all equal to $0$.
\end{itemize}

Let $\T^{Sp_{2n}}$ be the ring of invariants of $Sp_{2n}$ in 
$\T$.
The following is then straightforward from \fullref{thm:thm4.1} and 
\cite[Proposition 21]{G}.

\begin{thm}
\label{thm:thm5.1}
$\T^{Sp_{2n}}$ is the subring of
$\F_p[\ka_{n,1},\ldots,\ka_{n,n-1}]$ given by:
\begin{enumerate}
\item
for $p=2$, $\T^{Sp_{2n}}=\F_p[\ka_{n,0},\ldots,\ka_{n,n-1}]$;
\item
for $p>s2$,
\begin{enumerate}
\item
as a vector space over $\F_p$, $\T^{Sp_{2n}}$
has a basis $\{\ka^R_n | R \in \R'\}$;
\item
as a module over polynomial algebra
$\F_p[\ka_{n,0}^2,\ka^p_{n,1},\ldots,\ka^p_{n,n-1}]$,
$\T^{Sp_{2n}}$ is freely generated by $\{1,\ka^R_n | R \in 
\R''\}$.
\end{enumerate}
\end{enumerate}
\end{thm}

\bibliographystyle{gtart}
\bibliography{link}

\end{document}